 
\input amstex
\input epsf
\loadmsbm
\documentstyle{gen-j}
\NoBlackBoxes
\font\tit=cmr10 scaled\magstep3
\font\sfa=cmss10 scaled\magstep1

\def\({[}
\def\){]}
\font\bsl=cmbxsl10
\define\bslu{\hbox{\bsl u}}
 \catcode`\@=11
\def\bqed{\ifhmode\unskip\nobreak\fi\quad
  \ifmmode\blacksquare\else$\m@th\blacksquare$\fi}
\def\mex{\qopname@{mex}}

\catcode`\@=\active
\UseAMSsymbols
\vglue1.0cm
\topmatter
\title{}\endtitle

\endtopmatter

\centerline{{\tit Infinite cyclic impartial games}}\bigskip\bigskip
\centerline{{\sfa Aviezri S. Fraenkel
\footnote{\tt fraenkel\@wisdom.weizmann.ac.il\qquad
http://www.wisdom.weizmann.ac.il/\~{}fraenkel}
and Ofer Rahat\footnote{\tt ofer\@wisdom.weizmann.ac.il}}}\bigskip
\centerline{Department of Applied Mathematics and Computer Science} 
\centerline{Weizmann Institute of Science} 
\centerline{Rehovot 76100, Israel}
\vskip1.2cm
{\bf Abstract.} We define the family of {\it locally path-bounded\/} 
digraphs, which 
is a class of infinite digraphs, and show that on this class 
it is relatively easy to compute an optimal strategy (winning or 
nonlosing); and realize a win, when possible, in a finite number of 
moves. This is done by proving that the Generalized Sprague-Grundy 
function exists uniquely and has finite values on this class. 

\vskip0.8cm
\document
\centerline{\bf 1. Introduction}\medskip
We are concerned with {\it combinatorial games\/}, which, for our 
purposes here, comprise $2$-player games with perfect information, no 
chance moves and outcome restricted to (lose, win), (draw, draw) for 
the two players. A draw position is a position in the game such that 
no win is possible from it, but there exists a next move which 
guarantees, for the player making it, not to lose. You {\sl win\/} a 
game by making a last move in it. A game is {\it impartial\/} if for 
every position 
in it, both players have the same set of next moves; otherwise it's 
{\it partizan}. Nim is impartial, chess partizan. A game is {\it cyclic\/} 
if it contains cycles (the possibility of returning to the same position), 
or loops (pass-positions). These notions, slightly changed here, can be 
found in \cite{BCG82}. It is clear that a necessary (yet not sufficient) 
condition for the existence of draw positions is that the game be cyclic. 

For Partizan cyclic games, see \cite{Con78},\ \cite{Sha79},\ 
\cite{Fla81},\ \cite{FrTa82},\ \cite{Fla83}; 
finite impartial cyclic games are discussed only briefly in \cite{BCG82},  
\cite{Con76}. Particular finite impartial cyclic games are analysed in 
\cite{FrTa75}, \cite{FrKo87}. Infinite impartial games are treated briefly at 
the end of \cite{Smi66}, where both the ``generalized Sprague-Grundy
function'' $\gamma$, defined below, and its associated ``counter 
function'' were permitted to be transfinite ordinals. 

Our purpose here is to define a certain class of infinite digraphs on 
which $\gamma$ assumes only finite values, but the counter function may 
contain transfinite ordinal values. The motivation for doing this is 
based, in part, on the following 
considerations. It is easier to compute with finite than with transfinite 
ordinals. Often the structure of the digraph is such that the 
$\gamma$-function itself suffices to provide a winning strategy, 
without the need of an additional counter function (\S4). The 
$\gamma$-function always provides at least a nonlosing strategy; 
it's for consummating a win that the counter function may be needed. 
For consummating a win the ``generalized Nim-sum'' of a finite set 
of $\gamma$-values is required (\S4). The generalized Nim-sum is based 
on the binary expansion of ordinals. It's easy to see that every 
ordinal, finite or transfinite, has a unique expansion as a 
{\sl finite\/} sum of powers of ordinals (based on the greedy algorithm 
and the fact that the ordinals are well-ordered --- see 
\cite{Sie58 \rm{XIV, \S19}}). 
For example, $\omega=2^{\omega}$. We do not wish to enter here into 
the question of the computational complexity of computing with 
transfinite ordinals. But it seems possible that it's easier to 
compare the size of ordinals with each other, which suffices for 
counter function values, than to compute and work with their binary 
expansions, as needed for the $\gamma$-values.  

The connection between games and digraphs is simple: with any impartial 
game $\Gamma$ we associate a digraph $G=(V,E)$ where $V$ is the set 
of positions of $\Gamma$ and $(a,b)\in E$ if and only if there is a move 
from position $a$ to position $b$. It is called the {\it game-graph\/} 
of $\Gamma$. We identify games with their corresponding game-graphs, 
game positions with digraph vertices and game moves with digraph edges, 
using them interchangeably. It is thus natural to define a {\it cyclic\/} 
digraph as a digraph, finite or infinite, which may contain cycles or 
loops. 

In \S2 we provide basic tools needed for the statement and proof of
the result (Theorem~1), and \S3 contains the proof. An example 
demonstrating Theorem~1 is given in the final \S4.\medskip

\centerline{\bf 2. Preliminaries}\medskip

The subset of nonnegative integers is denoted by $\Bbb Z^0$, and the 
subset of positive integers by $\Bbb Z^+$. 

Given a digraph $G=(V,E)$. For any vertex $u\in V$, the set of 
{\it followers\/} of $u$ is $F(u)=\{v\in V:(u,v)\in E\}$. A 
vertex $u$ with $F(u)=\emptyset$ is a {\it leaf}. The set of 
{\it predecessors\/} of $u$ is $F^{-1}(u)=\{w\in V:(w,u)\in E\}$. 
A {\it walk\/} in $G$ is any sequence of vertices $u_1,u_2,\ldots$, 
not necessarily distinct, such that $(u_i,u_{i+1})\in E$, i.e., 
$u_{i+1}\in F(u_i)$ ($i\in\Bbb Z^+$). Edges may be repeated. A {\it path} 
is a walk with all vertices distinct. In particular, there's no 
repeated edge in a path. The {\it length\/} of a path is the number 
of its edges. If every path in $G$ has finite length, then $G$ is 
called {\it path-finite\/}. If there exists $b\in\Bbb Z^0$ such that 
every path in $G$ has length $\le b$, then $G$ is {\it path-bounded}.  

\definition{Definition 1} A cyclic digraph is {\it locally path-bounded\/} 
if for every vertex $u_i$ there is a bound $b_i(u_i)=b_i\in\Bbb Z^0$ such 
that the length of every (directed) path emanating from $u_i$ doesn't 
exceed  $b_i$. The integer $b_i$ is the {\it local path bound\/} of $u_i$. 
\enddefinition

Note that every path-bounded digraph is locally path-bounded, and 
every locally path-bounded digraph is path-finite. But neither of 
the two inverse relationships needs to hold. Our main result is 
concerned with locally path-bounded digraphs. 

Given a digraph $G=(V,E)$. The {\it Generalized Sprague-Grundy
function\/}, also called $\gamma$-{\it function\/}, is a mapping 
$\gamma\colon V\rightarrow\Bbb Z^0\cup\{\infty\}$, where the symbol $\infty$
indicates a value larger than any natural number. If $\gamma(u)=\infty$, 
we say that $\gamma(u)$ is infinite. We wish to define
$\gamma$ also on certain subsets of vertices. Specifically: $\gamma\bigl(
F(u)\bigr)=\bigl\{\gamma(v)<\infty\colon v\in F(u)\bigr\}$. If $\gamma(u)
=\infty$ and $\gamma\bigl(F(u)\bigr)=K$, we also write $\gamma(u)=\infty(
K)$. Next we define equality of $\gamma(u)$ and $\gamma(v)$: if $\gamma(u
)=k$ and $\gamma(v)=\ell$ then $\gamma(u)=\gamma(v)$ if one of the
following holds: (a)~$k=\ell<\infty$; (b)~$k=\infty(K)$, $\ell=\infty
(L)$ and $K=L$. We also use the notations
$$V^f=\bigl\{u\in V\colon\gamma(u)<\infty\bigr\},\qquad V^{\infty}=
V\setminus V^f,$$ 
where for any finite subset $S\subset\Bbb Z^0$, the {\it M}inimum 
{\it EX}cluded value $\mex$ is defined by 
$$\mex S=\min (\Bbb Z^0\setminus S)=\text{minimum term in}\ \Bbb Z^0\ 
\text{not in}\ S.$$
It is also convenient to introduce the notation 
$$\gamma'(u)=\mex\gamma\bigl(F(u)\bigr)=\mex\{\gamma(v)<\infty:v\in F(u)\}.
\tag 1$$

We need some device to tell the winner where to go when we use the 
$\gamma$-function. For example, suppose that there is a token on vertex $u$ 
(Fig.~1). It turns out that it's best for the player moving now to go 
to a position with $\gamma$-value $0$. There are two such values: one 
(the leaf) is an immediate win, and the other $(v)$ is only a nonlosing 
move. This digraph may be embedded in a large digraph where it's not 
clear which option leads to a win. The device which overcomes this 
problem is a counter function, as used in the following definition. 
For realizing an optimal strategy, we will normally select a follower 
of least counter function value with specified $\gamma$-value. The 
counter function also enables us to prove assertions by induction. 

\definition{Definition 2} Given a cyclic digraph $G=(V,E)$. A function
$\gamma\colon V\rightarrow\Bbb Z^0\cup\{\infty\}$ is a $\gamma$-{\it
function\/} with {\it counter function\/} $c\colon V^f\rightarrow J$,
where $J$ is any infinite well-ordered set, if the following three
conditions hold:

{\bf A}. If $\gamma(u)<\infty$, then $\gamma(u)=\gamma^{\prime}(u)$.\par
{\bf B}. If there exists $v\in F(u)$ with $\gamma(v)>\gamma(u)$, then 
there exists
$w\in F(v)$ satisfying $\gamma(w)=\gamma(u)$ and $c(w)<c(u)$.\par
{\bf C}. If $\gamma(u)=\infty$, then there is $v\in F(u)$
with $\gamma(v)=\infty(K)$ such that $\gamma^{\prime}(u)\notin K$.
\enddefinition

{\bf Remarks.} 
\roster
\item"{$\bullet$}" In {\bf B} we have necessarily $u\in V^f$; and we may have 
$\gamma(v)=\infty$ as in {\bf C}. 
\item"{$\bullet$}" To make condition {\bf C} more accessible, we state 
it also in the following equivalent form:

{\bf C'}. If for every $v\in F(u)$ with $\gamma(v)=\infty$ there is 
$w\in F(v)$ with $\gamma(w)=\gamma^{\prime}(u)$, then $\gamma(u)<\infty$.
\item"{$\bullet$}" If condition {\bf C'} is satisfied, then 
$ \gamma(u)<\infty$, and so by {\bf A}, $\gamma(w)=\gamma^{\prime}(u)=
\gamma(u)$. 
\item"{$\bullet$}" To keep the notation simple, we write 
$\infty(0)$, $\infty (1)$, $\infty (0,1)$ etc., for $\infty (\{ 0\} )$,
$\infty(\{ 1\} )$, $\infty (\{ 0,1\} )$, etc.\endroster

The $\gamma$-function was first defined in \cite{Smi66}. It was found 
independently in \cite{FrPe75}. The simplified version given above, 
and two other versions, appear in \cite{FrYe86}. Since this function 
is not well-known, we repeated its definition above. The $\gamma$-function 
exists uniquely on any finite cyclic digraph, but its associated counter 
function exists nonuniquely. Here and below, when we discuss the 
existence of $\gamma$, we mean its existence as a {\sl finite\/} 
ordinal (besides the special value $\infty$).

We are now ready to state our main result. 

\proclaim{\bf Theorem 1} Every locally path-bounded digraph\/ $G=(V,E)$ 
has a unique\/ $\gamma$-function with an associated counter function; 
and for every\/ $u\in V^f$, $\gamma(u)$ doesn't exceed the length of a 
longest path emanating from\/ $u$. 
\endproclaim

\centerline{\bf 3. The Proof}\medskip

We wish to examine some properties of path-finite and locally 
path-bounded digraphs. To begin with, is it clear that for a path-finite 
digraph, if $v\in F(u)$, then every path emanating from $v$ is not 
longer than any path emanating from $u$?

Perhaps it is clear, but it's also wrong: 
in a path-finite graph, every path originating at some vertex $u$ and
continuing to its ultimate end, terminates at a vertex $v$, where $v$ is
either a leaf or a predecessor of some $w$ on the path. Thus a path of
minimum length emanating from $u$ in Fig.~1 terminates at the leaf,
whereas a path of maximum length beginning at $u$ terminates at $y$. 
It has length 3. But a maximal-length path emanating from $v\in F(u)$ 
clearly has length 4. 
\midinsert\vskip 7pt
\centerline{\epsfxsize 4.0cm \epsffile{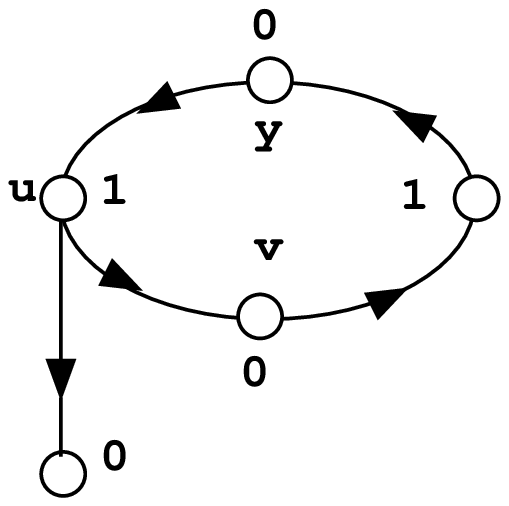}}
\botcaption{Figure 1. {\rm The numbers are $\gamma$-values.}}\endcaption
\endinsert\par
However, having embarked on a path $u_0,u_1,\dots,u_n$ of maximum length 
$n$, the maximum path length from any vertex $u_i$ encountered on 
it is $n-i$ $(0\le i\le n)$. 

If a digraph 
$G=(V,E)$, possibly with infinite paths, has no leaf, then the label 
$\infty$ on all the vertices is evidently a $\gamma$-function: {\bf A} 
and {\bf B} are satisfied vacuously, and {\bf C} is satisfied with 
$\gamma'(u)=0$ for all $u\in V$. If $G$ has a leaf, then some of the 
vertices have a $\gamma$-function, such as the leaf and its predecessors, 
but possibly $\gamma$ doesn't exist on some of the vertices. For 
the case where $\gamma$ exists on a subset $V'\subseteq V$, we define 
$\gamma'(u)=\mex\{\gamma(F(u)) : F(u)\subseteq V'\}$. Since $F(u)$
may, nevertheless, be infinite for any vertex $u$ in a locally path-bounded 
digraph, it is not clear a priori that $\gamma'(u)$ exists. The 
following lemma takes care of this point.

\proclaim{\bf Lemma 1} Let\/ $u$ be any vertex with local path bound\/ 
$b$ in a locally path-bounded digraph\/ $G=(V,E)$. Then\/ $\gamma'(u)$ 
exists\/ $($i.e., it is a nonnegative integer\/$)$, and in fact,\/ 
$\gamma'(u)\leq b$.\endproclaim  

\demo{\bf Proof} We consider two cases.

(i) Suppose that $u$ has finite $\gamma$-value $m$. Then $\gamma'(u)$ 
exists, and in fact, $\gamma'(u)=\gamma(u)=m$ by {\bf A}. Moreover, 
there exists $u_1\in F(u)$ with $\gamma(u_1)=m-1$, there exists 
$u_2\in F(u_1)$ with $\gamma(u_2)=m-2,\ldots,$ there exists $u_m\in F(u_{m-1})$
with $\gamma(u_m)=0$. Then $u,u_1,\ldots,u_m$ is a path of length $m$, so
$m\leq b$. (The path may continue beyond $u_m$, but in any case 
$\gamma'(u)=m\leq b$.)

(ii) Suppose that $u$ has either no $\gamma$-value or value $\infty$. It 
suffices to show that if $v\in F(u)\cap V^f$, then $\gamma(v)<b$. Indeed, 
$|F(u)|$ may be infinite, and $F(u)$ may contain vertices with no 
$\gamma$-value. But if $\gamma(v)<b$ for all $v\in F(u)\cap V^f$, then 
clearly $\gamma'(u)$ exists and $\gamma'(u)\le b$. Note that 
we cannot use the argument of case~(i) directly on $v$, since a path from 
$v$ may be longer than a path from $u$, as we just saw. So suppose there is 
$v_0\in F(u)\cap V^f$ with $\gamma(v_0)=n\geq b$. As in case (i), there is a
path $v_0,v_1,\ldots,v_n$ of length $n$ with $\gamma(v_i)=n-i\ 
(i\in \{0,\ldots,n\})$. Then $u,v_0,v_1,\ldots,v_n$ is a walk of length $n+1>b$
emanating from $u$. Hence it cannot be a path. But $v_i\ne v_j$ for all 
$i\ne j$, since $\gamma(v_i)\ne\gamma(v_j)$. Hence $v_j=u$ for some 
$j\in\{0,\ldots,n\}$. The contradiction is that $v_j$ does and $u$ 
doesn't have a finite $\gamma$-value. Thus $\gamma(v)<b$ for all 
$v\in F(u)\cap V^f$, hence $\gamma'(u)$ exists, and in fact, 
$\gamma'(u)\leq b$.\bqed\enddemo

\demo{\bf Proof of Theorem 1} Let $V'\subseteq V$ be a maximal subset 
of vertices on which $\gamma$ exists, together with an associated 
counter function $c$, subject to the following additions to {\bf B} and 
{\bf C} of Definition~2: 
$$\multline
\text{If}\ \gamma(u)<\infty\ \text{and there is}\ v\in F(u)\cap V_{\nu},\ \\
\text{then there is}\ w\in F(v)\ \text{with}\ \gamma(w)=\gamma(u),\ c(w)<c(u),
\endmultline\tag 2$$
$$\multline
\text{if}\ \gamma(u)=\infty,\ \text{then there is}\ v\in F(u)\ 
\text{with}\ \gamma(v)=\infty\ \\ 
\text{such that}\ w\in F(v)\cap V_{\nu}\implies \gamma'(w)\ne\gamma'(u),
\endmultline\tag 3$$
where $V_{\nu}=V\setminus V'$. (In (3) we have $\gamma'(w)\ne\gamma'(u)$, 
instead of $w\in F(v)$ and $\gamma(w)\ne\gamma'(u)$ in {\bf C}.) In 
addition we require:
$$\text{If}\ \gamma(u)=\infty\ \text{with}\ \gamma'(u)=l,\ 
\text{then}\ \gamma'(v)\ge l\ \text{for all}\ v\in V_{\nu}.\tag 4$$ 
The subset $V'$ is maximal in the sense that adjoining any $u\in V_{\nu}$ 
into $V'$ violates either Definition~2, or (2) or (3) or (4). If 
$V_{\nu}\ne\emptyset$, let $u\in V_{\nu}$. By Lemma~1, 
$\gamma'(u)=k$ exists for some $k\in\Bbb Z^0$. It follows that there 
is a minimum value 
$m=\min\{k\in\Bbb Z^0 : u\in V_{\nu},\ \gamma'(u)=k\}$. Let 
$K=\{u\in V_{\nu}:\gamma'(u)=m\}$. Then 
$V_{\nu}\ne\emptyset\implies K\ne\emptyset$. We consider four cases. 

{\bf Case 1.} For every $u\in K$ we have $m\in\gamma'(F(u))$, where, 
consistent with (1), 
$$\multline
\gamma'(F(u))=\{\gamma'(v):v\in F(u)\}=\{\mex \gamma(F(v)):v\in F(u)\}\\ 
=\{\mex\{\gamma(w)<\infty:w\in F(v), v\in F(u)\}\}.\endmultline$$
Note that $u\in K$, $v\in F(u)\implies\gamma(v)\ne m$ by the definition 
of $\mex$, so $v\in F(u)$, 
$\gamma'(v)=m\implies v\in V_{\nu}$, in fact, $v\in K$. Thus putting 
$\gamma(u)=\infty$ for all $u\in K$ satisfies {\bf C}, and is also 
consistent with (3); and with (4) by the minimality of $m$. Furthermore, 
it doesn't violate {\bf A}, and is consistent with {\bf B} by (2). 
This contradicts the maximality of $V'$. 

We may thus assume henceforth that there exists $u\in K$ such that 
$$m\notin\gamma'(F(u)).\tag 5$$ 

{\bf Case 2.} There exist $u\in K$ and $v\in F(u)$ with 
$\gamma(v)=\infty$, such that for every $w\in F(v)$, either 
$\gamma(w)\ne m$, or $w\in V_{\nu}$ with $\gamma'(w)\ne m$. 
Putting $\gamma(u)=\infty$ is clearly consistent with {\bf C}, and (3); 
and it doesn't violate {\bf A}. In view of (2), also {\bf B} is satisfied. 
This contradicts the maximality of $V'$. So we may assume that
$$\multline  
\forall u\in K\ \text{and}\ \forall v\in F(u)\ 
\text{with}\ \gamma(v)=\infty, \\ 
\exists w\in F(v)\ (\text{with}\ 
\gamma(w)=m\ \text{or}\ w\in V_{\nu}\ \text{with}\ \gamma'(w)=m).
\endmultline$$

We subdivide this into the following two cases:
$$\multline
\exists u\in K\ \text{such that}\ \forall v\in F(u)\ 
\text{with}\ \gamma(v)=\infty, \\ 
\exists w\in F(v)\ \text{with}\ \gamma(w)=m,
\endmultline\tag 6$$
or
$$\multline  
\forall u\in K\ \text{and}\ \forall v\in F(u)\ 
\text{with}\ \gamma(v)=\infty, \\ 
\exists\ \text{no}\ w\in F(v)\cap V',\ \text{but}\ 
\exists w\in F(v)\cap V_{\nu}\ \text{with}\ \gamma'(w)=m.
\endmultline\tag 7$$

{\bf Case 3.} (6) holds. We repeat that for any $u\in K$, since 
$\gamma'(u)=m$, $u$ has no follower with $\gamma$-value $m$. Suppose 
that there exists $y\in F^{-1}(u)$ with $\gamma(y)=m$. Then by (2), 
there exists $v\in F(u)$ with $\gamma(v)=m$, contradicting $\gamma'(u)=m$. 
Thus putting $\gamma(u)=m$ is consistent with {\bf A}. It is also consistent 
with (3): putting $\gamma(u)=m$ could presumably increase $\gamma'(y)$
for some $y\in F^{-1}(u)\cap V_{\nu}$, and thus upset (3) for the 
value $\gamma(z)=\infty$ of some grandparent $z=F^{-1}(F^{-1}(y))$ 
of $y$. Now by (4), $\gamma'(u)\ge\gamma'(z)$. If indeed $\gamma'(y)$ 
increased, then for the new value we have $\gamma'(y)>\gamma'(u)$, 
so $\gamma'(y)>\gamma'(z)$ and $\gamma(z)=\infty$ remains unaffected. 
Consistency with {\bf C} thus follows from (3) which becomes 
{\bf C} when $u$ is labeled $m$. Since $y\in F^{-1}(u)\implies\gamma(y)\ne m$, 
as we saw at the beginning of this case, the potential adverse effect 
on any grandparent $z$ of $y$ considered above, cannot happen. 

We now show that also {\bf B} holds. Suppose first that $F(u)\subseteq V'$. 
For every $v\in F(u)$ for which $\gamma(v)>m$, there exists $w\in F(v)$ 
with $\gamma(w)=m$. This follows from {\bf A} if $\gamma(v)<\infty$, 
and from (6) if $\gamma(v)=\infty$. It remains to define $c(u)$ 
sufficiently large so that $c(w)<c(u)$. This will be done below. 

In view of the minimality of $m$ and by (5), the second possibility 
is that for every $v\in F(u)$ for which $v\in V_{\nu}$, we have 
$\gamma'(v)>m$. For every such $v$ there exists $w\in F(v)$ with 
$\gamma(w)=m$ by the definition of $\mex$. Again we have to define 
$c(u)$ sufficiently large to satisfy $c(w)<c(u)$. 

Let $S=\{v\in F(u) : \gamma(v)>m\}\cup \{v\in F(u) : v\in V_{\nu},\ 
\gamma'(v)>m\}$. 
We have just seen that for every $v\in S$ there is $w\in F(v)$ with 
$\gamma(w)=m$. Put $T=\{w\in F(v) : v\in S,\ \gamma(w)=m\}$. Let 
$c(u)$ be the smallest ordinal $>c(w)$ for all $w\in T$. Then also 
{\bf B} is satisfied. This contradicts the presumed maximality of $V'$.

Note that the case $F(u)\subseteq V\setminus V^{\infty}$ satisfies (6) 
vacuously, and so is also included in the present case.  
\medskip 

{\bf Case 4.} (7) holds. If (7) holds nonvacuously, then as in Case~1, 
putting $\gamma(u)=\infty$ for all $u\in K$ is consistent with {\bf C}, 
(3) and the other conditions. This contradicts again the maximality of 
$V'$. Hence $K=\emptyset$, and so also $V_{\nu}=\emptyset$.\medskip 

Whenever $\gamma$ exists on a digraph, finite or infinite, it exists 
there uniquely. See \cite{Fra$\ge$99}. Finally, if $b$ is the 
local bound of $u\in V$, then $\gamma'(u)\le b$ by Lemma~2. Hence if 
$u\in V^f$, then $\gamma(u)\le b$ by {\bf A}.\bqed\enddemo\medskip

\centerline{\bf 4. An Example}\medskip

We specify below a locally path-bounded digraph $G=(V,E)$ on some of 
whose vertices we place a finite number of tokens. A move consists of 
selecting a token and moving it to a follower. Multiple occupancy of 
vertices is permitted. The player first unable to move loses, and the 
opponent wins. If there is no last move, the outcome is a draw. 

For any $r\in\Bbb Z^0$, a {\it Nim-heap\/} of size $r$ is a digraph 
with vertices $u_0,\dots,u_r$ and edges $(u_j,u_i)$ for all 
$0\le i<j\le r$. In depicting $G$ (Fig.~2), we use the convention 
that bold lines and the vertices they connect constitute a Nim-heap, 
of which only adjacent (bold) edges are shown, to avoid cluttering 
the drawing. Thin lines denote ordinary edges. 

All the horizontal lines are thin, and each vertex $u_i$ on this 
horizontal line connects via a vertical thin edge to a Nim-heap 
$G_i$ pointing downwards, of size $\lfloor(4i+8)/3\rfloor ,\ i\in\Bbb Z^0$.
From $u_i$ there also emanates a Nim-heap $H_j$ of size 
$j=\lfloor (i+2)/3\rfloor$ pointing upward. Each $G_i$ has a back edge 
to its top vertex forming a cycle of length $i+1$. Thus $G_0$ has a 
loop at the top of its Nim-heap (of size $2$). There is an additional 
back edge to the vertex $u_i$ on the horizontal line, forming a cycle 
of length $i+3$. 

\midinsert\vskip 7pt
\centerline{\epsfysize 13.0cm \epsffile{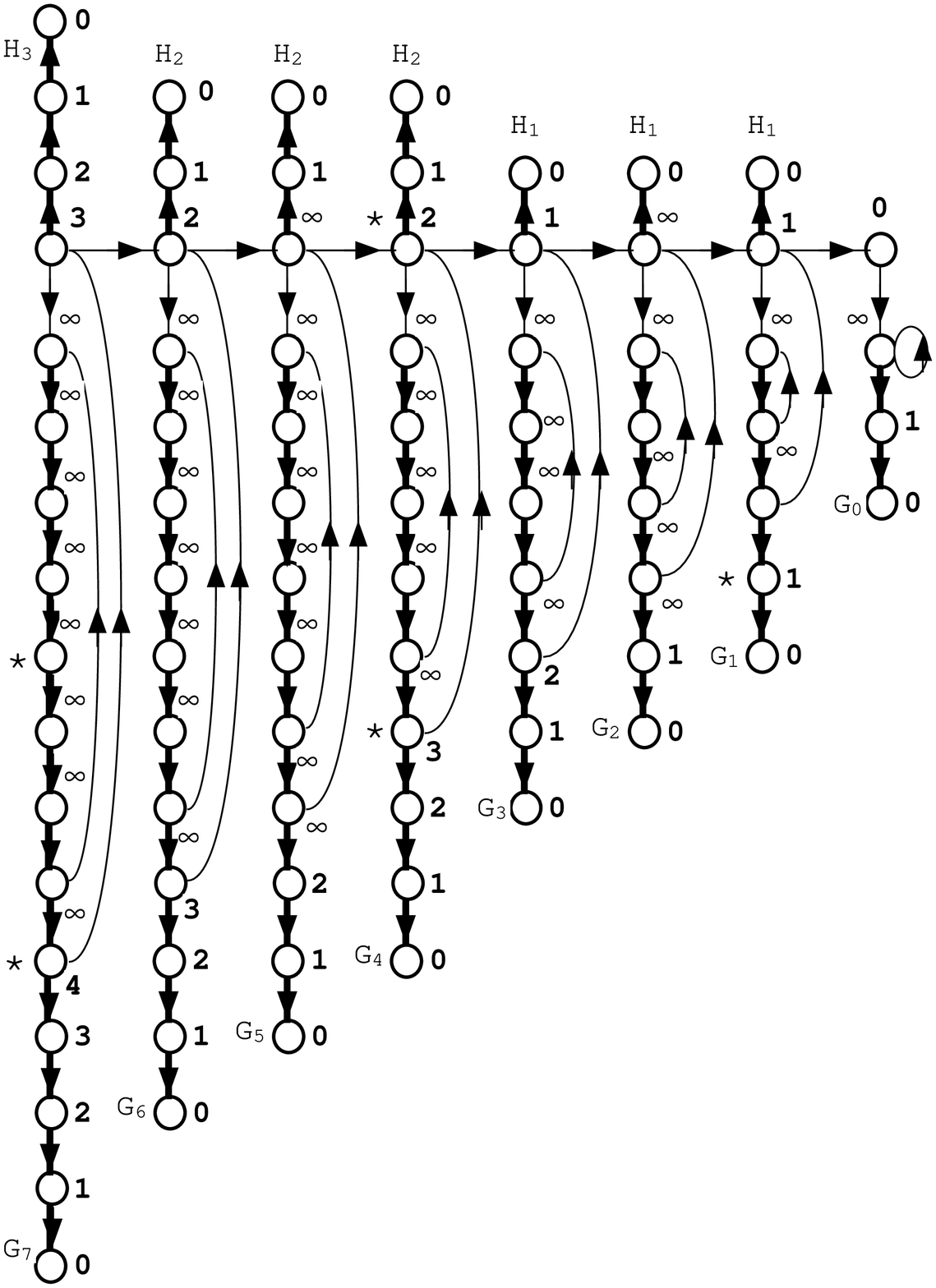}}
\botcaption{Figure 2. {\rm The tail-end of a locally path-bounded 
digraph.}}\endcaption
\endinsert\par
From any vertex $u$ on the horizontal line there is a longest path, 
via $G_i$, of length $\lfloor(4i+11)/3\rfloor$, and the other vertices 
have shorter maximal length. Thus $G$ is locally path-bounded. But it 
is not path-bounded, since $i$ can be arbitrarily large.\medskip

{\bf Sample Problem.} Compute an optimal strategy for the 5-token 
game placed on the $5$ starred vertices of $G$.\medskip 

To solve this problem, we introduce the generalized Nim-sum 
(\cite{Smi66}, \cite{FrYe86}, \cite{Fra$\ge$99}). 
For any nonnegative integer $h$ we write $h=\sum_{i\ge 0}h^i2^i$ 
for the binary encoding of $h$ ($h^i\in\{0,1\})$. If $a$ and $b$ 
are nonnegative integers, then their {\it Nim-sum} $a\oplus b=c$, 
also called {\it exclusive or, XOR}, or {\it addition over\/} GF$(2)$, 
is defined by
$c^i\equiv a^i+b^i\ (\hbox{mod}\ 2),\ c^i\in\{0,1\}\quad(i\ge 0).$

The Generalized Nim-sum of a nonnegative integer $a$ and $\infty
(L)$, for any finite subset $L\subset\Bbb Z^0$, is defined by
$a\oplus\infty (L)=\infty (L)\oplus a=\infty (L\oplus a),$ where 
$L\oplus a=\{\ell \oplus a\colon \ell \in L\}$. The Generalized Nim-sum of
$\infty (L_1)$ and $\infty (L_2)$, for any finite subsets $L_1$, $L_2$, 
is defined by
$\infty (L_1)\oplus\infty (L_2)=\infty (L_2)\oplus\infty(L_1)=\infty(
\emptyset).$
Clearly the Generalized Nim-sum is associative and $a\oplus a=0$ for 
every $a$. 

Given any finite or infinite game $\Gamma$, we say informally that 
a $P$-{\it position\/} is any position $u$ from which the {\it P}revious 
player can force a win, that is, the opponent of the player 
moving from $u$. An $N$-{\it position\/} is any position $v$ from which 
the {\it N}ext player can force a win, that is, the player who moves 
from $v$. A $D$-position is any position $u$ from which neither player 
can force a win, but has a nonlosing next move. The set of all $P$-, $N$- 
and $D$-positions is denoted by $\Cal P,\ $ $\Cal N$ and $\Cal D$ 
respectively. 

For any finite multiset $\bslu=(u_1,\dots,u_n)$ of vertices of $G$ on which 
tokens reside, one token on each $u_i$, we then have the result 
\cite{Fra$\ge$99}: 

\proclaim{\bf Proposition} 
The\/ $P$-, $N$- and $D$-labels of\/ $\bslu=(u_1,\dots,u_n)$ in any 
locally path-bounded digraph $G$ are given by
$$\eqalign{\Cal P&=\{\bslu\in V\colon\sigma (\bslu )=0\} ,\quad\Cal D =\{
\bslu\in V\colon\sigma (\bslu )=\infty (K),~~0\notin K\}\cr
\Cal N&=\{\bslu\in V\colon 0<\sigma(\bslu)<\infty\}\cup\{\bslu\in V\colon
\sigma(\bslu)=\infty (K),~~0\in K\} .\bqed}$$\endproclaim\medskip

We are now ready to solve the above problem, by observing that the symbols 
appearing on Fig.~2 are the $\gamma$-values of $G$. Simply check that 
they satisfy the conditions of Definition~2. In particular, {\bf B} 
of Definition~2 is satisfied if every vertex on the horizontal line 
with $\gamma$-value $<\infty$ gets a counter-value between $\omega$ 
and $\omega 2$, and every vertex in the Nim-heaps with $\gamma$-value 
$<\infty$ is assigned a counter value $<\omega$, which is clearly feasible. 

For the $5$ starred vertices we then have $1\oplus 3\oplus 2\oplus 4
\oplus\infty(0,1,2,3,4)=4\oplus\infty(0,1,2,3,4)=\infty(4,5,6,7,0)$, 
which contains $0$, hence the position is in $\Cal N$. Thus the player 
moving from this position can win by going to a position of Nim-sum 
$0$, namely, pushing the token on the infinity label to $4$. Indeed the 
resulting Nim-sum is $1\oplus 3\oplus 2\oplus 4\oplus 4=0$. 

We remark that any tokens on two vertices with $\gamma$-value $\infty$ 
is a draw position, no matter where the other tokens are, if any. Also 
note that for realizing a win in this game we do not really need a 
counter function.\medskip

\centerline{\bf Epilogue}\medskip

We have defined locally path-bounded digraphs, and shown that the 
generalized Sprague-Grundy function $\gamma$ exists on such digraphs 
with finite, though not necessarily bounded, values. Of course local 
path-boundedness is only a sufficient condition for the existence of 
$\gamma$. Any finite or infinite digraph without a leaf, satisfies 
trivially $\gamma(u)=\infty$ for all its vertices $u$. 

A large part of combinatorial game theory is concerned, however, 
with digraphs which do have leaves. If we exclude digraphs without 
leaves, then Theorem~1 is, in a sense, best possible. 

Consider the digraph $G$ which consists of a vertex $u$, and 
$F(u)=\{u_0,u_1,\dots\}$, where, for all $i\in\Bbb Z^0$, $u_i$ 
is the top vertex of a Nim-heap of size $i$, so $\gamma(u_i)=i$. 
Any path emanating from $u$ has the form $u,u_i,\dots$ for some 
$i$; its length is $i+1$. Paths not emanating from $u$ are shorter. 
Thus $G$ is path-finite. But $\gamma(F(u))=\{0,1,\dots\}$, so 
$\gamma(u)$ cannot assume any finite value.\medskip

\centerline{\bf References}\medskip

1. [BCG82] E. R. Berlekamp, J. H. Conway and R. K. Guy \(1982\), {\it Winning
Ways\/} (two volumes), Academic Press, London.

2. [Con76] J. H. Conway \(1976\), {\it On Numbers and Games\/}, Academic Press,
London.

3. [Con78] J. H. Conway \(1978\), Loopy Games, {\it Ann. Discrete Math.} 
{\bf 3}: Proc. Symp. Advances in Graph Theory, Cambridge Combinatorial
Conf. (B. Bollob\'as, ed.), Cambridge, May 1977, pp. 55--74.

4. [Fla81] J. A. Flanigan \(1981\), Selective sums of loopy partizan 
graph games, {\it Internat. J. Game Theory\/} {\bf 10}, 1--10.

5. [Fla83] J. A. Flanigan \(1983\), Slow joins of loopy games, {\it J. Combin.
Theory\/} (Ser.~A) {\bf 34}, 46--59.

6. [Fra$\ge$99] A. S. Fraenkel [$\ge$1999], {\it Adventures in Games and 
Computational Complexity}, to appear in Graduate Studies in Mathematics,
Amer. Math. Soc., Providence, RI.

7. [FrKo87] A. S. Fraenkel and A. Kotzig \(1987\), Partizan octal games: 
partizan subtraction games, {\it Internat. J. Game Theory\/} {\bf 16}, 
145--154.

8. [FrPe75] A. S. Fraenkel and Y. Perl \(1975\), Constructions in combinatorial
games with cycles, {\it Coll. Math. Soc. J\'anos Bolyai\/}, {\bf
10}:~Proc. Internat. Colloq. on Infinite and Finite Sets, Vol.~2 (A.
Hajnal, R. Rado and V. T. S\'os, eds.) Keszthely, Hungary, 1973,
North-Holland, pp. 667--699.

9. [FrTa75] A. S. Fraenkel and U. Tassa \(1975\), Strategy for a class of 
games with dynamic ties, {\it Comput. Math. Appl.} {\bf 1}, 237--254.

10. [FrTa82] A. S. Fraenkel and U. Tassa \(1982\), Strategies for compounds of
partizan games, {\it Math. Proc. Camb. Phil. Soc.} {\bf 92},
193--204.

11. [FrYe86] A. S. Fraenkel and Y. Yesha \(1986\), The generalized 
Spra\-gue-Grundy function and its invariance under certain mappings, 
{\it J. Combin. Theory\/} (Ser.~A) {\bf 43}, 165--177.

12. [Sha79] A. S. Shaki \(1979\), Algebraic solutions of partizan games with 
cycles, {\it Math. Proc. Camb. Phil. Soc.} {\bf 85}, 227--246.

13. [Sie58] W. Sierpi\'nski [1958], {\it Cardinal and Ordinal Numbers\/}, 
Hafner, New York.

14. [Smi66] C. A. B.\ Smith \(1966\), Graphs and composite games, 
{\it J. Combin. Theory\/} {\bf 1}, 51--81. Reprinted in slightly modified 
form in: {\it A Seminar on Graph Theory\/} (F. Harary, ed.), Holt, 
Rinehart and Winston, New York, NY, 1967.

\enddocument